\documentclass[a4paper]{amsart}
\usepackage{stmaryrd}
\usepackage{xr}
\usepackage{amscd, amssymb, calc, ifthen, mathrsfs}
\usepackage[matrix,arrow,curve,cmtip]{xy}
\usepackage{amsxtra}
\usepackage{eucal}
\usepackage{color}
\usepackage[normalem]{ulem}
\usepackage{amsxtra}
\usepackage{parskip}
\usepackage{mathabx}
\usepackage{extarrows}

\newcommand{\Fr}{\mathrm{Fr}}

\let \: \colon
\newcommand{\uni}{\mathrm{u}}
\newcommand{\Euni}{\mathscr{E}}
\newcommand{\xiuni}{\xi^\uni}
\newcommand{\xhat}{\wedge}
\newcommand{\ordn}{{\mathrm{o}}}
\newcommand{\mathbold}{\mathbb}
\newcommand{\bF}{{\mathbold F}}

\newcommand{\bZ}{{\mathbold Z}}

\newcommand{\bN}{{\mathbold N}}

\DeclareMathOperator{\Spec}{Spec}
\DeclareMathOperator{\Spf}{Spf}

\newcommand{\cn}{\mathrm{can}}

\newcommand{\longmap}{{\,\longrightarrow\,}}
\newcommand{\longequals}{{\,\xlongequal{\phantom{\rightarrow}}\,}}
\newcommand{\longlabelmap}[1]{{\,\buildrel #1\over\longrightarrow\,}}
\def\longisomap{{\,\buildrel \sim\over\longrightarrow\,}} 

\DeclareMathOperator{\colim}{\mathrm{colim}}
\newcommand{\pr}{{\mathrm{pr}}}

\newcommand{\sO}{{\mathcal{O}}}

\newcommand{\Hom}{\mathrm{Hom}}

\newcommand{\Aff}{\mathsf{Aff}}
\newcommand{\Afftil}{\Aff\sptilde}

\newcommand{\id}{{\mathrm{id}}}

\def\isomap{{\buildrel \sim\over\rightarrow}} 

\newcommand{\rightlabelxyarrows}[2]{{\ar@<0.7ex>^-{#1}[r]\ar@<-0.7ex>_-{#2}[r]}}

\newcommand{\displaylabelfork}[6]{{	\entrymodifiers={+!!<0pt,\fontdimen22\textfont2>}
	\def\objectstyle{\displaystyle}
\xymatrix{{#1} \ar^-{#2}[r] & {#3} \ar@<0.7ex>^-{#4}[r]\ar@<-0.7ex>_-{#5}[r] & {#6}}}}

\newcommand{\predisplaylabelfork}[6]{{{#1} \ar^-{#2}[r] & {#3} \ar@<0.7ex>^-{#4}[r]\ar@<-0.7ex>_-{#5}[r] & {#6}}}

\newtheoremstyle{mythm}{}{}%
  {\itshape}
  {}
  {\bfseries}
  {}
  { }
  {\thmnumber{#2.\hspace{1.5mm}}\thmname{#1}\thmnote{ {\mdseries(#3)}}.}

\newtheorem*{theorem*}{Theorem} 

\newtheoremstyle{intro}{}{}%
  {\itshape}
  {}
  {\bfseries}
  {}
  { }
  {\thmname{#1}\thmnumber{ #2}\thmnote{ #3}.}

\numberwithin{equation}{subsection}
\theoremstyle{mythm}
\newtheorem{theorem}[subsection]{Theorem}

\theoremstyle{intro}

\address{Mathematical Sciences Institute\\ Australian National University\\Canberra ACT 0200\\Australia}

\email{james.borger@anu.edu.au}

\address{Korteweg-de Vries Instituut\\ Universiteit van Amsterdam\\1090 GE Amsterdam\\The Netherlands}
\email{l.r.gurney@uva.nl}

\title{Canonical lifts of families of elliptic curves}

\author[J.~Borger, L.~Gurney]{James Borger, Lance Gurney}
\thanks{This work was supported the Australian Research Council.}

\begin{document}

\begin{abstract}

We show that the canonical-lift construction for ordinary elliptic curves over perfect fields of characteristic
$p>0$ extends uniquely to arbitrary families of ordinary elliptic curves, even over $p$-adic formal schemes. In
particular, the universal ordinary elliptic curve has a canonical lift. The existence statement is largely a
formal consequence of the universal property of Witt vectors applied to the moduli space of ordinary elliptic
curves, at least with enough level structure. As an application, we show how this point of view allows for more
formal proofs of recent results of Finotti and Erdo\u{g}an.

\end{abstract}

\maketitle

\section{}

Fix a prime number $p$. Let $W$ denote the usual, $p$-typical Witt vector functor. Let $R$ be a ring in which
$p$ is nilpotent, and write $S=\Spec R$. Let $W_n(S)$ denote $\Spec W_n(R)$, and let $W(S)$ denote the
direct limit $\colim_n W_n(S)$. We take this limit in the category of sheaves of sets on the category of affine
schemes with respect to the \'etale topology. One could say that $W(S)$ is the correct version $\Spec W(R)$,
a construction which, as we discuss below, does not have good properties.

We say an elliptic curve $E$ over $S$ is ordinary when all fibers of $E$, necessarily over points of residue
characteristic $p$, are ordinary. For any morphism $f\:S'\to S$, we will write $E_{S'}$, or $f^*(E)$, for 
the base change $S'\times_S E$ regarded as an elliptic curve over $S'$ in the evident way. 

The purpose of this paper is to prove the following:

\begin{theorem*}

There is a unique way of lifting ordinary elliptic curves $E$ over affine schemes $S$ on which $p$ is
nilpotent to elliptic curves $\widetilde{E}$ over $W(S)$ which is compatible with base change in $S$ and has the
property that $\widetilde{E}$ admits a Frobenius lift $\psi\:\widetilde{E}\to F^*(\widetilde{E})$,
where $F$ is the usual Witt vector Frobenius map $F\:W(S)\to W(S)$.

\end{theorem*}

Note that the requirement here that $S$ is affine is only to simplify the exposition. We will remove it below
and allow $S$ to be any $p$-adic formal scheme, or even what we call a $p$-adic sheaf. See 
section~\ref{sec-thm-statement} for the final statement of the theorem and further details.

We call $\widetilde{E}$ the canonical lift of $E$. In the case $S=\Spec k$ where $k$ is a perfect field of 
characteristic $p$, our canonical lift agrees with the usual one by the remarks in~\ref{subsec:perfect-base}.

\section{Background on sheaves}

In this section and the next, we define $W(S)$, the infinite-length Witt vector construction when $S$ is a
scheme, and even when $S$ is more general. The reason there is something to do is that while $\Spec W_n(R)$ is
a well-behaved construction, the naive infinite-length analogue $\Spec W(R)$ is not. For instance, some basic
geometric facts like the theorem in section~\ref{sec-W-for-schemes} are not true for the naive infinite-length
construction. Thus to handle Witt vectors of infinite length geometrically, one needs some way of retaining the
information of the projective system of the $W_n(R)$ instead of crudely passing to the limit $W(R)=\lim_n
W_n(R)$ in the category of rings.

One way of doing this is to view $W(R)$ as a topological ring with the inverse-limit topology, each $W_n(R)$
being discrete. In the special case where $R$ is a perfect $\bF_p$-algebra, this topology is an adic topology
with ideal of definition generated by $p$. We can therefore consider its formal spectrum $\Spf W(R)$, and in
this way, the theory of formal schemes can accommodate a satisfactory theory of Witt vector constructions of
infinite length, as long as $R$ is perfect. For general rings $R$, the inverse-limit topology on
$W(R)$ is not an adic topology, and so the theory of formal schemes cannot even get started.

Another way of proceeding, also standard,
is to view $\Spec W_n(R)$ as a sheaf of sets on the category of affine schemes, and to define $W(\Spec R)$ to be
the direct limit $\colim_n \Spec W_n(R)$ in this category. Since all limit and colimit constructions are as well
behaved in categories of sheaves as they are in the category of sets, there is good reason to be confident in
this approach. Indeed it works without problem and is the one we will take.

There is a second issue in defining $W(S)$, which is how general we should allow $S$ to be. It will be
convenient later if we can allow $S$ to be the moduli space of ordinary elliptic curves (with sufficient level
structure) over rings on which $p$ is nilpotent. But this does not exist in the category of schemes. It does
exist in the category of formal schemes, where one can construct it as the ordinary locus in the $p$-adic
completion of the moduli space of all elliptic curves. However it also exists in the category of sheaves
mentioned above. Since we will be using this category anyway, it will be simpler to take that approach and
forget about formal schemes entirely. It also has the side benefit of working for objects $S$ much
more general than $p$-adic formal schemes.

\subsection{Objects representable over sheaves}
\label{subsec-objects-over-sheaves}

Let $\Aff$ denote the category of affine schemes, and let $\Afftil$ denote the category of sheaves of sets on
$\Aff$ with respect to the \'etale topology. (See SGA4~\cite{SGA4.2}, exp.\ VII.) Any scheme $S$ can be viewed
as an object of $\Afftil$ via the functor is represents $\Spec R\mapsto \Hom(\Spec R,S)$. This is a fully
faithful embedding and we will regard the category of schemes as a full subcategory of $\Afftil$ in this way
without further comment.

For any $S\in\Afftil$, let $\Afftil_S$ denote the category of sheaves $X$ equipped with a map $X\to S$, where
the morphisms are morphisms over $S$. If $T\to S$ is a morphism, and $X\in\Afftil_S$, let $X_{T}$ denote the
sheaf $T\times_S X$ together with the morphism $\pr_T\:T\times_S X\to T$. We can then define familiar
scheme-theoretic structures on objects $X$ over $S$ by using affine test schemes. For example, an elliptic
curve over $S$ is a sheaf $X$ over $S$ together with for every morphism $T\to S$ from an affine scheme $T$, an
elliptic curve structure on $X_T$ (that is, an isomorphism to $X_T$ from the sheaf represented by a usual
elliptic curve over $T$) with the property that for any morphism $T'\to T$ of affine schemes the
isomorphism $(X_T)_{T'} \isomap X_{T'}$ is a morphism of usual elliptic curves over $T'$. Note that it also
follows that the sheaf $X$ over $S$ admits a unique group structure (over $S$) inducing the group structures on
the elliptic curves $X_T$ over $T$.

Descent theory for $\Afftil$ goes through for general sheaf-theoretic reasons. Suppose $S'\to S$ is an
epimorphism in $\Afftil$; in other words, every section of $S$ lifts \'etale locally to one of $S'$. Then any
object of $\Afftil_{S'}$ equipped with descent data to $S$ descends to a unique (up to unique isomorphism)
object of $\Afftil$ over $S$. This remains true if we are interested in objects $X$ with additional structure,
as long as that structure is of an \'etale-local nature on affine schemes. This is the case for ordinary
elliptic curves, as they satisfy effective descent for the \'etale topology, and is the only example we will
need.

There is a generalization of this which we will also use. Suppose we are given a presentation of an object
$S\in\Afftil$ as a colimit: $\colim_i S_i\isomap S$. Then the category of objects over $\colim_i S_i$ is
equivalent to the category of compatible families of objects over the $S_i$. Indeed, an object $X$ over $S$
gives rise to a compatible family of objects $X_i=S_i\times_S X$ over the $S_i$, where compatible means that
we are given morphisms $X_i\to X_j$ lying over each morphism $S_i\to S_j$ such that the induced maps $X_i\to
S_i\times_{S_j}X_j$ are isomorphisms. Conversely, to such a compatible family $X_i$ we associate the object
$X=\colim_i X_i$. As above, the equivalence between objects over a colimit and compatible families also holds
for objects with additional structure of an \'etale-local nature, such as ordinary elliptic curves. For general
reasons, any sheaf $S$ admits such a presentation where each $S_i$ is an affine scheme, and so an ordinary
elliptic curve (say) over $S$ is equivalent to a compatible family of ordinary elliptic curves over this
diagram of affine schemes. More importantly, one often considers objects $S$ which are given as the colimit of
some specific family of affine schemes, and then one can describe an elliptic curve over $S$ as a compatible
family of elliptic curves over this specific family. So an elliptic curve over an object given in these terms
really is an accessible object.

\subsection{$p$-adic sheaves}
\label{sec-p-adic-sheaves}
We will say that a sheaf $S\in\Afftil$ is \emph{$p$-adic} if it is isomorphic to a colimit $\colim_i S_i$ of
affine schemes $S_i$ on which $p$ is nilpotent. (We require nothing of the indexing diagram of the colimit
other than that it is small. In other language, this means that it is a set and not a proper class.) For
example, a scheme is $p$-adic if and only if $p$ is locally nilpotent on it. In fact, any $p$-adic sheaf $S$
has the property that $p$ is nilpotent on any affine scheme admitting morphism to $S$ (and the converse is true
up to issues of set-theoretic size just mentioned).

If $T$ is a scheme, put
	$$
	T_n=\Spec \bZ/p^{n+1}\bZ \times_{\Spec \bZ} T.
	$$
Then we call
	$
	\widehat{T}=\colim_n T_n
	$ 
the $p$-adic completion of $T$. It is a $p$-adic sheaf because each $T_n$ is a scheme over $\bZ/p^{n+1}\bZ$
and hence a colimit of affine schemes over $\bZ/p^{n+1}\bZ$.

We can do the same with any $p$-adic formal scheme $T$, and this defines a fully faithful embedding of the
category of $p$-adic formal schemes into the category of $p$-adic sheaves.

\section{Background on Witt vectors}

\subsection{Universal property of Witt vectors for rings} 
\label{subsec-can-lift-rings} 
Let $A$ be a
$p$-torsion free ring with an endomorphism $\psi\:A\to A$ such that $\psi(a)\equiv a^p \bmod pA$.
Let $R$ be any ring such that $W(R)$ is $p$-torsion free. (This does not hold in general, but it does hold
at two opposite extremes---when $R$ is a perfect $\bF_p$-algebra and when $R$ is $p$-torsion free.)
Then any ring map $g\:A\to R$ lifts to a unique ring map $\widetilde{g}$
	$$
	\xymatrix{
	A \ar^{g}[dr]\ar@{-->}^{\widetilde{g}}[rr]
		&& W(R) \ar[dl] \\
	& R
	}
	$$
which is Frobenius equivariant, meaning that
$F\circ \widetilde{g} = \widetilde{g}\circ \psi$.
In other words, the image of $\widetilde{g}(a)$ under the ghost map $W(R)\to R^\bN$  is
\begin{equation}
	\label{eq:ghost-comp}
	\langle g(a),\, g(\psi(a)),\, g(\psi^2(a)),\,\dots\rangle.
\end{equation}

A proof using the traditional definition of Witt vectors can be found in Lazard's 
book~\cite{Lazard:formal-groups-book}, p.\
215. We remark however that the traditional definition of Witt vectors will be irrelevant in this paper---it is
only the universal property that matters. So a preferable alternative would be to take a definition of Witt
vectors making the universal property obvious. For such a development, one can see Joyal~\cite{Joyal:Witt} for a
concise account or section 1 of~\cite{Borger:BGWV-I} for a more extensive one.

Observe that $\widetilde{g}\:A\to W(R)$ agrees with the composition $W(g)\circ\widetilde{\id_A}\:A\to W(A)\to W(R)$. We
can therefore define a canonical map $\widetilde{g}$ without any restrictions on $R$. Indeed, we simply define
$\widetilde{g}=W(g)\circ\widetilde{\id_A}$. One might call $\widetilde{g}$ the canonical lift of $g$. We emphasize that
while we have dropped all assumptions on $R$ here, we are still assuming $A$ is $p$-torsion free. We also
emphasize that while $\widetilde{g}$ is a canonically defined
Frobenius equivariant lift of $g$, without any assumptions on $R$ it might not be the unique one.

If we now write $Y=\Spec A$, then for any ring $R$, we have functorial maps $Y(R)\to Y(W(R))$ given by
$g\mapsto\widetilde{g}$. In other words, if $Y$ is a $p$-torsion-free moduli space parametrizing objects of some
given type, a Frobenius lift on $Y$ defines a way of canonically lifting objects defined over $R$ to objects
over $W(R)$. We might say that a class of objects has a theory of canonical lifts whenever their moduli space is
$p$-torsion free and has a Frobenius lift. (For a little more along these lines see section \ref{subsec:epilogue}.) Indeed, this is the principle we will apply below. But we will need
it in a slightly modified form because in our example, the Frobenius lift $\psi$ exists only on the $p$-adic
completion of the moduli space, which is not a scheme, but the $p$-adic completion of a scheme,
in the sense of section~\ref{sec-p-adic-sheaves}. The form of the universal property we need we be given in
section~\ref{subsect-can-lifts-sheaves}.

\subsection{Witt vectors for schemes}
\label{sec-W-for-schemes}
Write $W_n\:\Aff\to\Aff$ for the functor defined by $W_n(\Spec R)=\Spec W_n(R)$. 
The following theorem is fundamental:

\begin{theorem*}
	The functor $W_n\:\Aff\to\Aff$ 
	preserves \'etale maps, \'etale covering families, and \'etale base change.
	In particular, $W_n$ is continuous in the \'etale topology.
\end{theorem*}

The first general argument was given in van der Kallen~\cite{van-der-Kallen:Descent}, (2.4), but was written
only for the `big' Witt vector functor. For a proof for the $p$-typical Witt vector functor considered here,
one can see~\cite{Borger:BGWV-I}, section 9.2. One can also see the earlier paper by Langer and
Zink~\cite{Langer-Zink:dRW}, appendix A; note that while the results there are stated only under some
finite-type assumptions, the general case can be deduced by a limiting argument. For all our applications,
however, it is enough to consider Witt vectors of rings in which $p$ is nilpotent, and in this context the
results stated in Langer--Zink apply without any finiteness assumptions.

This allows us to extend $W_n$ to $\Afftil$, by SGA4~\cite{SGA4.1}, exp.\ III. Indeed, because $W_n$ is
continuous, for any sheaf $S$, the presheaf $U\mapsto S(W_n(U))$, for any $U\in \Aff$, is a sheaf. This defines
a functor $W_{n*}\:\Afftil\to\Afftil$, and it has a left adjoint $W_n^*$. Finally, $W_n^*$ extends $W_n$ from
$\Aff$ to $\Afftil$ in the sense that we have canonical isomorphisms $W_n^*(\Spec R)\isomap\Spec W_n(R)$. So
from now on, we will often abusively write $W_n=W_n^*$.

Since $W_n\:\Afftil\to\Afftil$ is a left adjoint, it preserves colimits. Therefore any presentation 
$f\:\colim_i \Spec R_i\isomap S$ gives rise to a presentation 
$$
W_n(f)\:\colim_i \Spec W_n(R_i) \longisomap W_n(S).
$$ 
Therefore by the remarks in section~\ref{subsec-objects-over-sheaves}, an
elliptic curve over $W_n(S)$ can be understood as a compatible system of elliptic
curves over the $W_n(S_i)$, and similarly for any other kind of object having a local nature.

Observe that if $S$ is a $p$-adic sheaf, then so is $W_n(S)$. Indeed, it is sufficient (and necessary) to
observe that if $p$ is nilpotent in a ring $R$, then it is nilpotent in $W_n(R)$. One way to show this
is to observe that the comonad structure map $W_{m+n}(\bF_p)\to
W_n(W_m(\bF_p))$ makes $W_n(\bZ/p^{m+1}\bZ)$ into a $\bZ/p^{m+n+1}\bZ$-algebra.

In fact, if $S$ is a scheme on which $p$ is locally nilpotent, then so is $W_n(S)$. We will not use this
below and mention it only for the reader who is more comfortable with schemes than abstract sheaf theory. As a
topological space, $W_n(S)$ agrees with $S$, and its structure sheaf $\sO_{W_n(S)}$ is given by the presheaf
$U\mapsto W_n(\sO_S(U))$. One way to show this is by using the fact that for any ring $R$ and any element $f\in
R$, we have $W_n(R[1/f])=W_n(R)[1/[f]]$, where $[f]$ denotes the Teichm\"uller lift of $f$. Alternatively, open
immersions are the same as \'etale monomorphisms and hence are preserved by $W_n$, by the theorem above.

Then we define 
	$$
	W(S)=\colim_{n} W_n(S).
	$$
We emphasize that this colimit is taken in $\Afftil$. Thus $W$ is the left adjoint of $\lim_n W_{n*}$. If
$S=\colim_i S_i$, then we have $W(S)=\colim_{n,i}W_n(S_i)$. Therefore if $S$ is a $p$-adic sheaf, then so is
$W(S)$. But $W(S)$ is typically not a scheme, even if $S$ is. (For instance, if $S=\Spec \bF_p$, then $W(S)$ is
the colimit of $\Spec \bZ/p^{n+1}\bZ$, which represents the functor sending $\Spec R$ to the one-point set if
$p$ is nilpotent in $R$ and to the empty set otherwise. This is representable by the formal scheme $\Spf \bZ_p$
but not by a scheme.) It is however still easy to work with. For example, if $T$ is an affine scheme, then the
set $\Hom(T,W(S))$ is the filtered colimit $\colim_{n}\Hom(T,W_n(S))$. (By for example SGA4~\cite{SGA4.2},
exp.\ VI, Thm.\ 1.23(ii) on p.\ 185.) So if $S$ is a scheme, for example, then any map $T\to W(S)$ factors
through the scheme $W_n(S)$ for some $n$.

Note that $W_n$ and $W$ preserve epimorphisms in $\Afftil$, as they are left adjoints. In particular,
if $S'\to S$ is an \'etale cover of schemes, then $W(S')\to W(S)$ is an epimorphism, and hence descent
is available.

\subsection{The universal property of Witt vectors for $p$-adic sheaves}
\label{subsect-can-lifts-sheaves}

Let $A$ be a $p$-torsion free ring with a Frobenius lift $\psi$. Let $Y_m=\Spec A/p^{m+1}A$, 
$\widehat{Y}=\colim_m Y_m$, as in section~\ref{sec-p-adic-sheaves}, and $\widehat{\psi}: \widehat{Y}\to \widehat{Y}$ be the Frobenius lift induced by $\psi$. Let $S$ be a $p$-adic sheaf and let $f\:S\to
\widehat{Y}$ be a morphism. Write $S=\colim_i \Spec R_i$, with $p$ nilpotent in each ring $R_i$. Then for each $i$,
there exists an $m_i$ such that the map $\Spec R_i\to \widehat{Y}$ factors through the inclusion $Y_{m_i}\to
\widehat{Y}$, thus inducing a map $A/p^{m_i+1}A\to R_i$. Let $g$ denote the composition $A\to A/p^{m_i+1}A \to
R_i$, and let $\widetilde{g}$ denote the canonical lift as defined in~\ref{subsec-can-lift-rings}. Then for each
$n$, the composition
	$$
	A\longlabelmap{\widetilde{g}} W(R_i) \longmap W_n(R_i)
	$$ 
factors through $A\to A/p^{N_{i,n}+1}A$ for some $N_{i,n}$ (in fact, for any $N_{i,n}\geq m_i+n$).
This defines a compatible family of maps 
	$$
	\Spec W_n(R_i)\longmap \Spec A/p^{N_{i,n}+1} \longequals Y_{N_{i,n}} \longmap \widehat{Y},
	$$
and hence, by the universal property of colimits, a map 
	$$
	\widetilde{f}\:W(S) \longequals \colim_{n,i}W_n(\Spec R_i)\longmap \widehat{Y},
	$$
which we will again call the canonical lift of $f\:S\to\widehat{Y}$. Note that as with the previously defined 
canonical lift maps, the map $\widetilde{f}$ is Frobenius equivariant, which is to say 
$\widetilde{f}\circ \widehat{\psi}=F\circ \widetilde{f}$.

\section{Statement of the theorem}
\label{sec-thm-statement}

The usual Witt vector Frobenius map $F\: W_{n+1}(R)\to W_n(R)$ induces functorial maps $F\:W_n(S)\to
W_{n+1}(S)$, for any $S\in\Afftil$, and upon taking colimits, maps $F\:W(S)\to W(S)$. It also satisfies the
relation $F(x)\equiv r(x)^p \bmod pW_n(R)$, where $r\:W_{n+1}(R)\to W_n(R)$ denotes the usual projection. So the
maps $F\:W(S)\to W(S)$ agree with the usual $p$-th power Frobenius map on the locus $\Spec\bF_p \times_{\Spec
\bZ}W(S)$.

Let us then say that a \emph{Frobenius lift} on an elliptic curve $E$ over $W(S)$ is a morphism $E\to F^*(E)$ of
elliptic curves over $W(S)$ restricting to the usual Frobenius map over the locus $\Spec \bF_p\times_{\Spec
\bZ}W(S)$ modulo $p$.

\begin{theorem*}
There is a unique way of lifting ordinary elliptic curves $E$ over $p$-adic sheaves $S$ 
to elliptic curves $\widetilde{E}$ over $W(S)$ such that
the construction $E\mapsto \widetilde{E}$ is compatible with base change in $S$ and such that
each $\widetilde{E}$ admits a Frobenius lift.
\end{theorem*}

We make some remarks to clarify the statement. First, the compatibility condition, which can be written
$(E_{S'})^\sim=\widetilde{E}_{W(S')}$ for any map $S'\to S$, is more properly expressed as a coherent family of
isomorphisms. Second, the uniqueness statement is to be understood as follows: if $E\mapsto\widecheck{E}$ is any
other such construction, then there is a unique family of isomorphisms $\widetilde{E}\to\widecheck{E}$, where $E$ runs
over all ordinary elliptic curves $E$ over all $p$-adic base sheaves $S$, which are Frobenius equivariant and
compatible with restriction of the base $S$. We emphasize that such a uniqueness statement does not apply to
lifts of a single elliptic curve or even all elliptic curves over a given base $S$ but only to the family of
all elliptic curves over all bases. However see section~\ref{thm:unique-Frob-lift} for a result in this
direction.

\section{Existence}
\label{sec-construction}
In this section, we construct the canonical lift functor. In the presence of enough level structure, it is
nothing more than the universal property of Witt vectors applied to the moduli space of ordinary elliptic curves
with its canonical Frobenius lift. In general, we will use a descent argument to pass from the case with
level structure to the general setting.

\subsection{Existence with level structure}
\label{elnf}

Let $Y(N)$ denote the moduli space of elliptic curves $E/S$ with full level-$N$ structure
$\xi\:(\bZ/N\bZ)^2\isomap E[N](S)$. We will assume $p\nmid N$ and that $N$ is large enough to make the moduli
problem representable. (So $N\geq 3$ is enough.) In this case, $Y(N)$ is a smooth affine scheme of relative
dimension one over $\Spec \bZ[1/N]$. Let $T$ denote the open subscheme of $Y(N)$ which is the complement of
the supersingular locus on the fiber over $p$, and let $Y(N)^\ordn$ denote its $p$-adic completion $\widehat{T}$,
in the sense of section~\ref{sec-p-adic-sheaves}. Then $Y(N)^\ordn$ is of the form $\colim_n
\Spec A_N/p^{n+1}A_N$, where $A_N$ is $p$-adically complete and $p$-torsion free. (Indeed, if we write $Y(N)=\Spec R$, then $A_N$ is the
$p$-adic completion of $R[Q^{-1}]$, for any subset $Q\subset R$ such that $\Spec R[Q^{-1}]/(p)$ is the
ordinary locus of the fiber of $Y(N)$ over $p$. It follows that $A_N$ is $p$-torsion free because $R$ smooth over
$\bZ[1/N]$ and hence $p$-torsion free.) Thus it is of the form needed to apply the universal property of
section~\ref{subsect-can-lifts-sheaves}

Recall the standard Frobenius lift $\psi$ on $Y(N)^\ordn$. Let $E$ be an ordinary elliptic curve over an affine
scheme $S$ on which $p$ is nilpotent, and let $E^\cn$ denote its canonical subgroup, the connected component of
its $p$-torsion sub-group-scheme (\cite{Katz:modular-forms}, Ch.\ 3). It is a finite flat closed
sub-group-scheme of $E$. Let $E/E^\cn$ denote the usual quotient (so not the quotient object in $\Afftil$, but
the quotient with respect to the fppf topology). Thus $E/E^\cn$ is also a family of ordinary elliptic curves
over $S$, and the quotient map $E\to E/E^\cn$ is faithfully flat. Further, if $\xi$ is a level-$N$ structure on
$E$, then its image $\bar{\xi}$ in $E/E^\cn$ is still a level-$N$ structure. We then let $\psi$ denote the map
$Y(N)^\ordn\to Y(N)^\ordn$ that, for any $S$, sends a an $S$-valued point $(E,\xi)$ to $(E/E^\cn,\bar{\xi})$.
It is a Frobenius lift because for ordinary elliptic curves over $\bF_p$-algebras, the connected component of
the $p$-torsion sub-group-scheme agrees with the kernel of Frobenius.

Now let $(E,\xi)$ be an ordinary elliptic curve with level-$N$ structure over a $p$-adic sheaf $S$.
This is the pull-back of the universal object $(\Euni,\xiuni)$ through a unique map
$c\:S\to Y(N)^\ordn$. Since the ring $A_N$ is $p$-torsion free, we can apply
the universal property of Witt vectors as given in section~\ref{subsect-can-lifts-sheaves} and write
	$$
	\widetilde{c}\: W(S) \to Y(N)^\ordn
	$$
for the morphism induced by $c$.
We then let $(E,\xi)^\sim$ denote the canonical lift of $(E,\xi)$, which is to say 
the pull-back $\widetilde{c}^*((\Euni,\xiuni))$ of the universal elliptic curve with level-$N$ structure.
The assignment 
	$$
	(E, \xi)\mapsto (E,\xi)^\sim
	$$
now defines our theory of canonical lifts for elliptic curves with full level-$N$ structure.
We could write $(E,\xi)^\sim = (\widetilde{E},\widetilde{\xi})$, but it will not be until the next section that
we know $\widetilde{E}$ is canonically independent of the choice of $\xi$.

\subsection{Existence in general}

We will use a standard descent argument. So let $Y(N,N)$ denote the moduli space of elliptic curves with a pair
of full level-$N$ structures $(\xi_1,\xi_2)$. Forgetting one or the other defines projections
$Y(N,N)\rightrightarrows Y(N)$, both which are finite \'etale. Now we proceed as we did with $Y(N)^\ordn$. Let
$Y(N,N)^\ordn$ denote the $p$-adic completion of the complement of the supersingular locus on the fiber over
$p$. Then $Y(N,N)^\ordn$ is of the form $\colim_n \Spec A_{N,N}/p^{n+1}A_{N,N}$, where $A_{N,N}$ is a $p$-adically complete and
$p$-torsion free ring, and it has a Frobenius lift, also denoted $\psi$, defined by sending an elliptic curve to
its quotient by the canonical subgroup, with the image level structures. This defines a theory of canonical
lifts for elliptic curves with pairs of full level-$N$ structure. It is compatible with the projections
$Y(N,N)\rightrightarrows Y(N)$ in the sense that taking the canonical lift commutes with forgetting each of the
level structures. This is simply because the Frobenius lift $\psi$ commutes with the projections.

We can now define the canonical lift of an arbitrary ordinary elliptic curve $E$ over a $p$-adic sheaf $S$. It
will in fact be enough to do this for affine schemes $S$, as the construction will be compatible with base
change along morphisms of affine schemes, so that if $S$ is an arbitrary $p$-adic sheaf and $E$ is an ordinary
elliptic curve over $S$, writing $S=\colim_i S_i$ as a colimit of affine schemes, we may define $\widetilde{E}$ over
$W(S)$ to be $\colim_i (E_{S_i})^\sim$ over $W(S)=\colim_i W(S_i)$, as in
section~\ref{subsec-objects-over-sheaves}.

So let $E$ be an ordinary elliptic curve over an affine $p$-adic scheme $S$. Let $S'$ be the universal cover of 
$S$ over which $E$ admits a level-$N$ structure $\xi$. The covering morphism
$S'\to S$ is then finite and \'etale, and we have a diagram
	$$
	\xymatrix{ S'\times_S S' \ar^{c_2}[r]\ar@<0.7ex>[d]\ar@<-0.7ex>[d] & Y(N,N)^\ordn
	\ar@<0.7ex>[d]\ar@<-0.7ex>[d] \\ S' \ar^{c_1}[r] & Y(N)^\ordn. } 
	$$ 
Both columns have the structure of a groupoid object in the category $\Afftil$.
The groupoid structure on the left column is unique because $S'\times_S S'$ is
an equivalence relation; and on the right column, it is the one with composition given by 
	$$
	(E,\xi_1,\xi_2)\circ (E',\xi'_1,\xi'_2) = (E,\xi_1,\beta^*(\xi'_2)).
	$$
whenever there is an isomorphism (necessarily unique) $\beta\:(E,\xi_2)\isomap (E',\xi'_1)$.
With respect to these groupoid structures, this diagram becomes a morphism of groupoid objects.
Because the
right column has Frobenius lifts which are compatible with the projections, the universal property of
Witt vectors then gives us a diagram
	$$ 
	\xymatrix{ W(S'\times_S S') \ar^-{\widetilde{c}_2}[r]\ar@<0.7ex>[d]\ar@<-0.7ex>[d] &
	Y(N,N)^\ordn \ar@<0.7ex>[d]\ar@<-0.7ex>[d] \\ W(S') \ar^-{\widetilde{c}_1}[r] & Y(N)^\ordn.} 
	$$ 
By the theorem in section~\ref{sec-W-for-schemes},
we have
	\begin{align*}
	W(S'\times_S S')&=\colim_n W_n(S'\times_S S') = \colim_n W_n(S')\times_{W_n(S)} W_n(S') \\
		&= W(S')\times_{W(S)}W(S'),
	\end{align*}
and so the diagram above can be identified with 
	$$
	\xymatrix{ 
	W(S')\times_{W(S)}W(S')
	\ar[r]\ar@<0.7ex>[d]\ar@<-0.7ex>[d] & Y(N,N)^\ordn \ar@<0.7ex>[d]\ar@<-0.7ex>[d] \\
	W(S') \ar[r] & Y(N)^\ordn, 
	} 
	$$ 
which is easily seen to be a morphism of groupoid objects again.
In other words, there is a family over $W(S')$ equipped with
descent data to $W(S)$. Since $W$ preserves epimorphisms, as explained above,
the map $W(S')\to W(S)$ is an effective descent morphism for elliptic curves. 
So we can define $\widetilde{E}$ to be the descended object over $W(S)$. It is well defined up to unique 
isomorphism, in the usual sense. As remarked earlier, it is at this point clear that the construction $E/S\mapsto \widetilde{E}/W(S)$ for affine schemes $S$ is compatible with base change and therefore can be extended to any ordinary elliptic curve over any $p$-adic sheaf $S$ and we do so without further comment.

Our construction of the canonical lift $\widetilde{E}$ appears to depend on the auxiliary choice of the level $N$.
One could easily show at this point that it does not, up to canonical isomorphism, but this is a consequence of
the uniqueness statement proved in section~\ref{sec-uniqueness}, and there is no need to establish it earlier.

\subsection{Remark: a stack-theoretic interpretation}\label{subsec:stack-interpetation}
The language of stacks is well-suited for expressing the descent argument above. For simplicity, we will
explain it for Witt vectors of finite length.

Given a stack $X$, let $W_{n*}(X)$ denote the fibered category sending any affine scheme $S$ to the category
$X(W_n(S))$. This is often called the arithmetic jet space or the Greenberg transform of $X$. It is
straightforward to show that $W_{n*}$ sends affine schemes to affine schemes, that it is a right adjoint and
hence sends groupoid objects to groupoid objects, and that it preserves \'etale morphisms.
(See~\cite{Borger:BGWV-II}.) Finally, by the theorem in~\ref{sec-W-for-schemes} the functor $W_n$ is continuous
in the \'etale topology. It then follows for general reasons that $W_{n*}$ takes stacks to stacks.
(See~\cite[Tag 04WC]{stacks-project}.)

The compatible Frobenius lifts on $Y(N)^\ordn$ and $Y(N,N)^\ordn$ can then be packaged as
a morphism of \'etale groupoid objects 
	$$
	\xymatrix{
	Y(N,N)^\ordn \ar[r] \ar@<0.7ex>[d]\ar@<-0.7ex>[d]
		& W_{n*}(Y(N,N)^\ordn) \ar@<0.7ex>[d]\ar@<-0.7ex>[d] \\
	Y(N)^\ordn \ar[r]
		& W_{n*}(Y(N)^\ordn).
	}
	$$
It therefore prolongs canonically to a morphism of the quotient stacks
	$$
	\xymatrix{
	Y(N,N)^\ordn \ar[r] \ar@<0.7ex>[d]\ar@<-0.7ex>[d]
		& W_{n*}(Y(N,N)^\ordn) \ar@<0.7ex>[d]\ar@<-0.7ex>[d] \\
	Y(N)^\ordn \ar[r] \ar[d]
		& W_{n*}(Y(N)^\ordn) \ar[d] \\
	Y(1)^\ordn \ar@{-->}[r]
		& W_{n*}(Y(1)^\ordn).
	}
	$$
This morphism of stacks is then nothing more than the family of canonical lift functors,
truncated at length $n$.
The image of an elliptic curve $E$ over $S$ is the elliptic curve $\widetilde{E}_{W_n(S)}$ over $W_n(S)$.

It is worth noting that Buium's work on differential modular forms~\cite{Buium:differential-modular-forms} also
touches on the connection between canonical lifts, moduli spaces, and arithmetic jet spaces.

\subsection{Frobenius lifts}
\label{subsec:Frob-lifts}

In this section, we define isomorphisms
	$$
	\eta_E\: \widetilde{E}/\widetilde{E}^\cn \longisomap F^*(\widetilde{E}),
	$$
where $F$ denotes the Frobenius map $W(S)\to W(S)$, such that the composition 
$$
\widetilde{E}\longrightarrow \widetilde{E}/\widetilde{E}^\cn\longisomap F^*(\widetilde{E})
$$
with the quotient map is a Frobenius lift, in the sense of section~\ref{sec-thm-statement}. First observe
that such an isomorphism $\eta_E$, if it exists, is necessarily unique. 
Indeed, to show two isomorphisms with the property above agree, it is enough to show they agree after base change to any
$W_n(S_i)$, where $S_i$ is any affine scheme mapping to $S$ on which $p$ is nilpotent;
but then the difference of two such maps would be a map of elliptic curves
which is zero over $\Spec \bF_p\times_{\Spec \bZ} W_n(S_i)$,
and hence zero over all of $W_n(S_i)$ by the rigidity theorem 
(\cite{Katz-Mazur}, thm.\ 2.4.2, p.\ 76).
Therefore local existence on the base will imply global existence.
In particular, we may assume that $E$ admits a level-$N$ structure and then, by base change, that 
$S=Y(N)^\ordn$ and that $E$ is the universal curve $\Euni$. 

The classifying morphism $c=\widetilde{\id}\:W(Y(N)^\ordn)\to Y(N)^\ordn$ for $\widetilde{\Euni}$
is Frobenius equivariant, by construction. In other words, there is a unique identification 
$c^*\psi^*(\Euni)=\psi^*F^*(\Euni)$ compatible with the level structure. 
We also have $\psi^*(\Euni)=\Euni/\Euni^\cn$, by the definition of $\psi$.
Thus we have identifications 
	$$
	\widetilde{\Euni}/\widetilde{\Euni}^\cn = c^*(\Euni/\Euni^\cn) = c^*\psi^*(\Euni)
	= F^*c^*(\Euni) = F^*(\widetilde{\Euni})
	$$
which are compatible with the level structure (and hence unique). Finally, the composition
$\widetilde{\Euni}\to\widetilde{\Euni}/\widetilde{\Euni}^\cn= F^*(\widetilde{\Euni})$ reduces to the relative $p$-th power
Frobenius map modulo $p$ because, writing $\widetilde{\Euni}_0=\Spec \bF_p \times_{\Spec \bZ} \widetilde{\Euni}$, the
actual Frobenius map $\widetilde{\Euni}_0\to \Fr^*(\widetilde{\Euni}_0)$ has kernel $\widetilde{\Euni}^\cn_0$ and is
compatible with the level structure.

From the stack-theoretic point of view, the map $\eta_E$ can be viewed as an invertible natural transformation
$\psi\circ c\to c\circ F$ and hence as providing a Frobenius equivariant structure on the morphism $c\:W(S)\to
Y(1)^\ordn$.

\subsection{Remark: Avoiding abstract sheaf theory}

It is possible to avoid abstract sheaf theory by working with compatible systems of elliptic curves. Indeed, as
explained in sections~\ref{subsec-objects-over-sheaves} and~\ref{sec-W-for-schemes}, an elliptic curve over
$W(S)$, where $S$ is a scheme on which $p$ is nilpotent, is equivalent to a compatible family of elliptic
curves over the $W_n(S)$, which are also schemes on which $p$ is nilpotent (and even of finite type over
$\bZ_p$ if $S$ is). Thus the canonical lift $\widetilde{E}$ can be viewed as a compatible system of elliptic curves
$\widetilde{E}_n$ over the schemes $W_n(S)$, which are not much harder to understand than $S$ itself.

Of course one could construct this compatible system directly, without going through Witt vectors of infinite
length. For each $n$, the truncated canonical lift $\widetilde{E}_n$ is constructed exactly as we did above with
Witt vectors of infinite length but using $W_n$ everywhere instead of $W$. One then shows directly that the
$\widetilde{E}_n$ form a compatible system, and we never have to leave the category of schemes.
The Frobenius morphism then becomes a family of morphisms 
	$$
	\widetilde{E}_n \longmap \widetilde{E}_n/\widetilde{E}_n^\cn \longisomap F^*(\widetilde{E}_{n+1})
	$$
of elliptic curves over the $W_n(S)$, where $F$ now denotes the truncated Witt vector Frobenius map 
$W_n(S)\to W_{n+1}(S)$.

\section{Uniqueness}
\label{sec-uniqueness}

To prove the uniqueness part of the theorem, we will need the following
result in the particular case of the universal elliptic curve. But since the result is just as easy to prove in
a more general form, we will do that.

\begin{theorem}\label{thm:unique-Frob-lift}
	Let $R$ be a $p$-adically complete ring such that $W(R)$ is $p$-torsion free. Let 
	$S=\colim_n \Spec R/p^{n+1}R$, and let $E$ be
	an ordinary elliptic curve over $S$. Suppose $X_1$ and $X_2$ are lifts of $E$ to $W(S)$ with
	Frobenius lifts $\psi_1$ and $\psi_2$. Then there is a unique Frobenius equivariant isomorphism 
	$\beta\:X_1\to X_2$ restricting to the identity on $E$.
\end{theorem}

We recall again that $W(R)$ is $p$-torsion free if $R$ is either $p$-torsion free or a perfect $\bF_p$-algebra.
Also note that by formal GAGA, elliptic curves over $S$ are equivalent to elliptic curves over $\Spec R$.

\begin{proof}
The Frobenius lifts
are morphisms $X_i \to F^*(X_i)$ of elliptic curves over $W(S)$ which reduce to the
Frobenius map modulo $p$, and hence are
(representable and) flat of degree $p$, by rigidity (\cite{Katz-Mazur}, thm.\ 2.4.2, p.\ 76).
Therefore the kernel is a (representable) finite flat
lift of the kernel of the Frobenius map and hence equals the canonical subgroup. Thus the Frobenius lifts
$\psi_i$ induce isomorphisms $\bar{\psi}_i\:X_i/X_i^\cn\isomap F^*(X_i)$.

We first consider the case where $E$ admits a level-$N$ structure $\xi$. Let $\xi_i$ denote the unique lift of
$\xi$ to $X_i$, and let $\alpha_i\:W(S)\to Y(N)^\ordn$ denote the classifying map for $(X_i,\xi_i)$. Now
observe that the isomorphisms $\bar{\psi}_i$ must preserve the level structure, where $X_i/X_i^\cn$ is given
the image level structure and $F^*(X_i)$ is given the pull-back level structure. Indeed, consider the elliptic
curves $\bar{X}_i=\Spec \bF_p\times_{\Spec \bZ} X_i$ over $\Spec \bF_p \times_{\Spec \bZ} W(S)$. Because
$p\nmid N$, the $N$-torsion is finite \'etale. Therefore the Frobenius map $\bar{X}_i \to \Fr_S^*(\bar{X}_i)$
preserves the level structure and, again because the $N$-torsion is finite \'etale, so does any lift of the
Frobenius map. In particular, the maps $\bar{\psi}_i$ preserve it.

It follows that the elliptic curves $X_i/X_i^\cn$
and $F^*(X_i)$ have the same classifying map $W(S)\to Y(N)^\ordn$. The classifying map for $X_i/X_i^\cn$ is
$\psi\circ\alpha_i$, by the definition of $\psi$, and that for $F^*(X_i)$ is $\alpha_i\circ F$.
Therefore, we have $\psi\circ\alpha_i=\alpha_i\circ F$, which is to say that the maps 
$\alpha_i:W(S)\to Y(N)^\ordn$ commute with the Frobenius lifts.

But at the same time, since both $(X_i,\xi_i)$ lift
$(E,\xi)$, the two compositions
	$$
	\xymatrix{S \ar[r] & W(S) \ar@<0.7ex>^{\alpha_1}[r]\ar@<-0.7ex>_{\alpha_2}[r] & Y(N)^\ordn}
	$$
agree. Therefore, writing $Y(N)^\ordn=\colim_n \Spec A/p^{n+1}A$ with $A$ $p$-adically complete,
we see that the two compositions
	$$
	\xymatrix{A \ar@<0.7ex>^-{\alpha_1^*}[r]\ar@<-0.7ex>_-{\alpha_2^*}[r] & W(R)\ar[r] & R}
	$$
agree. Further, the maps $\alpha_1^*$ and $\alpha_2^*$ are Frobenius equivariant because $\alpha_1$ and 
$\alpha_2$ are.
Since $A$ and $W(R)$ are $p$-torsion free, the universal property of Witt vectors implies $\alpha_1^*=\alpha_2^*$
and hence $\alpha_1=\alpha_2$.

When $E$ is arbitrary, the existence and uniqueness of $\beta$ follow formally because \'etale locally $\beta$
exists and is unique. We write out the details. Let $S'$ denote the universal \'etale cover of $S$ over which
$E$ admits a level-$N$ structure. Write $E'=S'\times_S E$ and $X'_i=W(S')\times_{W(S)} E_i$. Then each $X'_i$ is
a lift of $E'$ with a Frobenius lift $\psi'_i=F\times \psi_i$. Since each $X'_i$ also admits a level-$N$
structure, the construction above gives a canonical morphism $\beta'\:X'_1\to X'_2$ (and hence in fact a unique
one, because $N$ is large) commuting with the Frobenius lifts $\psi'_i$. Further, $\beta$ is equivariant with
respect to descent data. Indeed, write $S''=S'\times_S S'$ and $X''_i=W(S'')\times_{W(S)}X$. Then each $X''_i$
is a lift of $S''\times_S E$ with Frobenius lift $\psi''_i=F\times \psi_i$. Therefore there is a unique
Frobenius equivariant morphism $\beta''\:X''_1\to X''_2$. In particular, the two restrictions of $\beta'$ to
$W(S'')$ agree, which is to say that $\beta'$ is equivariant with respect to descent data.
\end{proof}

\subsection{Uniqueness of the canonical lift functor}
Suppose we have an assignment $E\mapsto \widecheck{E}$, sending elliptic curves $E/S$ to elliptic
curves $\widecheck{E}/W(S)$ with a Frobenius lift which is compatible with change of the base $S$. To give an 
isomorphism $\widecheck{E}\to\widetilde{E}$,
it is enough give isomorphisms locally on $S$ which are compatible with change of $S$. 
Therefore we may assume $E$ admits a level-$N$ structure. Let $d\:S\to Y(N)^\ordn$ denote the 
classifying morphism of $E$, and let $W(d)\:W(S)\to W(Y(N)^\ordn)$ denote the induced morphism.
Now since we have $Y(N)^\ordn=\colim_n \Spec A_N/p^{n+1}A_N$, where $A_N$ is $p$-adically complete and $p$-torsion free, we can apply the theorem above to the universal ordinary elliptic curve $\Euni$ over $Y(N)^\ordn$.
Therefore there is
a unique Frobenius equivariant isomorphism $\widecheck{\Euni}\to\widetilde{\Euni}$ restricting to the identity on 
$\Euni$. Since the constructions
$E\mapsto \widecheck{E}$ and $E\mapsto\widetilde{E}$ are compatible with restriction of $S$ (in the first case by 
assumption and in the second case by section~\ref{sec-construction}), the a family of morphisms
	$$
	\widecheck{E} \longequals W(d)^*(\widecheck{\Euni}) \longisomap W(d)^*(\widetilde{\Euni}) \longequals \widetilde{E},
	$$
is the unique family of Frobenius equivariant isomorphisms $\widecheck{E}\to\widetilde{E}$ 
which is compatible with restriction of $S$. 

In particular, the functor $E\mapsto\widetilde{E}$ constructed in section~\ref{sec-construction} is independent
of the choice of the level $N$, up to unique isomorphism.

\section{Remarks and applications}

\subsection{}
\label{subsec:perfect-base}

In the classical context where $S=\Spec k$ with $k$ a perfect field of characteristic $p$, it is known that the
Serre--Tate canonical lift has a Frobenius lift and is characterized, up to unique isomorphism, by this
property. (See Messing's book\cite{Messing:Barsotti-Tate}, p.\ 177, cor.\ (1.2) and p.\ 174, cor.\ (3.4).) The
theorem in section~\ref{thm:unique-Frob-lift} then implies that our canonical lifts are isomorphic to the
Serre--Tate canonical lifts by a unique Frobenius equivariant morphism.

Also observe that in this context, we only need a weak case of the theorem in~\ref{sec-W-for-schemes} in the
descent argument. This is because level structure exists over some finite extension $k'/k$, and then it is well
known that $W(k')/W(k)$ is finite \'etale and we have $W(k'\otimes_k k')=W(k')\otimes_{W(k)}W(k')$. Thus in the classical context, our argument really is little more than an application of the universal property of Witt
vectors.

\subsection{} While our approach has the benefit of allowing families in mixed characteristic, it also requires
them, even if one is only interested in canonical lifts of elliptic curves in characteristic $p$; and this
could be viewed as a drawback. The reason is simply that the special fiber $Y(N)^\ordn_{\bF_p}$ does not admit
a map from $W(Y(N)^\ordn_{\bF_p})$. Indeed, no nonempty Witt vector scheme maps to $\Spec \bF_p$.

\subsection{}

We can describe the $j$-invariant of the canonical lift of the 
universal family. Recall that $Y(N)^\ordn=\Spf A_N$ and $Y(N,N)^\ordn=\Spf
A_{N,N}$, where $A_N$ and $A_{N,N}$ are $p$-adically complete and $p$-torsion free, and that forgetting the
level structures induces a pair of Frobenius equivariant maps $A_N\rightrightarrows A_{N,N}$. Thus, writing
$A_1$ for the equalizer of these two maps, the map $({\id_{A_N}})^{\sim}\:A_N\to W(A_N)$ restricts to a map $s:
A_1\to W(A_1)$, and so the Frobenius lift on $A_N$ restricts to a Frobenius lift on $A_1$. Now write
$A_1=\bZ_p[j,1/f(j)]^\xhat$, where $j$ is an indeterminate identified with the $j$-function, $f(j)$ is a monic
polynomial whose roots lift the supersingular $j$-invariants, and $(\cdot)^\xhat$ denotes $p$-adic completion.
Then, as above, the image of $j$ under the map $s\:A_1\to W(A_1)$ has ghost components $\langle j,
\psi(j),\psi^{\circ 2}(j),\dots\rangle\in A_1^\infty$. It follows that if $E$ is any family of ordinary
elliptic curves, the ghost components of $j(\widetilde{E})$ are obtained by evaluating the universal expressions
$\psi^{\circ n}(j)\in\bZ_p[j,1/f(j)]^\xhat$ at $j=j(E)$.

One can similarly consider the usual Witt components
$(j_0,j_1,\dots)$ of $s(j)$. Then the Witt components of $j(\widetilde{E})$ are also obtained by evaluating
the universal expressions $j_n\in \bZ_p[j,1/f(j)]^\xhat$ at $j=j(E)$. For example, by the definition of the
ghost map, we have $\psi^{\circ n}(j)=\sum_{i=0}^n p^i j_i^{p^{n-i}}$, and so the Witt components
of the $j$-invariant of the canonical lift of an elliptic curve with $j$-invariant $j$ are given by 
	$$
	j_0=j, \quad j_1=(\psi(j)-j^p)/p, \quad j_2=(\psi(\psi(j))-j^{p^2})/p^2-(\psi(j)-j^p)^p/p^{p+1},
	$$
and so on. For elliptic curves over perfect fields of characteristic $p$, such a result was found by
Erdo\u{g}an~\cite{Erdogan:canonical-lifts} (strengthening earlier results of Finotti~\cite{Finotti}) by more
traditional means, but where the universal expressions are of course the reductions of ours modulo $p$. We note
that Erdo\u{g}an also constructs canonical lifts for families over perfect $\bF_p$-algebras.

\subsection{} 

We can do the same for elliptic curves over a general base.
Let $E$ be a family over
$\Spf R$ for a $p$-adically complete ring $R$. Then the $j$-invariant $j(\widetilde{E})$ is an element of 
$W(R)$, and as in~(\ref{eq:ghost-comp}) its ghost components are 
	$$
	\langle j(E_0), j(E_1), j(E_2),\dots \rangle =
	\langle j(E), \psi(j(E)), \psi^{\circ 2}(j(E)),\dots \rangle
	\in R^\infty,
	$$ 
where $E_0=E$ and $E_{n+1}=E_n/E_n^\cn$. This already determines $j({\widetilde{E}})$ if $R$ is $p$-torsion
free, since in that case the ghost map $W(R)\to R^\bN$ is injective. In general, one can lift $E$ arbitrarily to
a family $E'$ over a $p$-torsion-free extension $R'$ of $R$, and then $j(\widetilde{E})$ is the image under $
W(R')\to W(R)$ of $j((E')^\sim)$, which has the simple description above since $R'$ is $p$-torsion free.

This can be expressed in an especially pleasant way in
the particular case where $R$ is a perfect $\bF_p$-algebra $k$, and $R'$ is $W(k)$. Then one can show that the
projection $W(\pr_0)\:W(W(k))\to W(k)$ sends a Witt vector with ghost components $\langle
a_0,a_1,\dots\rangle\in W(k)^\infty$ to $\lim_{n\to\infty} F^{-n}(a_n)\in W(k)$, and so we have
	$$
	j(\widetilde{E})=\lim_{n\to\infty} F^{-n}(j(E'_n))=\lim_{n\to\infty} F^{-n}(\psi^{\circ n}(j(E'))),
	$$ 
where $E'$ is an arbitrary lift to $W(k)$ of $E$ and the $E'_n$ are as above.
If $k$ is a field with $p^r$ elements, then $F^r=\id_{W(k)}$ and
so we simply have $j(\widetilde{E})=\lim_{n\to\infty} j(E'_{rn})$.

\subsection{}

Amusingly, one can take canonical lifts of canonical lifts and so on, any number of times. Indeed, if $E$ is an
ordinary elliptic curve over a $p$-adic sheaf $S$, then $\widetilde{E}$ is also an ordinary elliptic curve over
the $p$-adic sheaf $W(S)$. However one can show that the double canonical lift $\widetilde{\widetilde{E}}$ is
canonically isomorphic to the pull-back of $\widetilde{E}$ through the composition, or plethysm, map $W(W(S))\to
W(S)$ of~\cite{Borger:BGWV-II}, (10.6.13), and similarly for the higher iterates. So these higher canonical
lifts are all determined by the first one in a way that has nothing to do with elliptic curves and only the
general theory of Witt vectors.

\subsection{Epilogue: On a general theory of canonical lifts}\label{subsec:epilogue}

We proposed above that a theory of canonical lifts for a class of objects parameterized by a moduli space $Y$
should be defined to be a Frobenius lift $\psi\:Y\to Y$, as long as $Y$ is $p$-torsion-free. In the absence of
this assumption, one should define a theory of canonical lifts to be a slightly stronger structure. Namely, the functor $W_*=\lim_n W_{n*}$ has a natural comonad structure
(coming from the maps $W_{m+n*}\to W_{m*}\circ W_{n*}$ of~\cite{Borger:BGWV-II}, (10.6.14)) and a theory of canonical lifts for the objects parametrized by $Y$ would then be a coaction of the comonad $W_*$ on $Y$. This is called is called a $\delta$-structure in Joyal~\cite{Joyal:Witt} and Buium~\cite{Buium:Diff-chars-over-p-adic}
and a $p$-typical $\Lambda$-structure in~\cite{Borger:BGWV-I}\cite{Borger:BGWV-II}.

This can also be done for other varieties of Witt vectors in the sense of~\cite{Borger:BGWV-II}, such as the
big Witt vectors. One would then
define a theory of canonical lifts, relative to the given variety of Witt vectors, to be an action of the
comonad $W_*$ on $Y$. In other words, a theory of canonical lifts for a class of objects should defined to be a $\Lambda$-structure on their moduli space. This was carried out by the second author \cite{Thesis-Gurney} in the case of elliptic curves with complex multiplication and the variety of big Witt vectors associated with an imaginary quadratic field.

\frenchspacing
\bibliography{references}

\begin{thebibliography}{10}

\bibitem{SGA4.1}
{\em Th\'eorie des topos et cohomologie \'etale des sch\'emas. {T}ome 1:
  {T}h\'eorie des topos}.
\newblock Springer-Verlag, Berlin, 1972.
\newblock S\'eminaire de G\'eom\'etrie Alg\'ebrique du Bois-Marie 1963--1964
  (SGA 4), Dirig\'e par M. Artin, A. Grothendieck, et J. L. Verdier. Avec la
  collaboration de N. Bourbaki, P. Deligne et B. Saint-Donat, Lecture Notes in
  Mathematics, Vol. 269.

\bibitem{SGA4.2}
{\em Th\'eorie des topos et cohomologie \'etale des sch\'emas. {T}ome 2}.
\newblock Springer-Verlag, Berlin, 1972.
\newblock S\'eminaire de G\'eom\'etrie Alg\'ebrique du Bois-Marie 1963--1964
  (SGA 4), Dirig\'e par M. Artin, A. Grothendieck et J. L. Verdier. Avec la
  collaboration de N. Bourbaki, P. Deligne et B. Saint-Donat, Lecture Notes in
  Mathematics, Vol. 270.

\bibitem{Borger:BGWV-I}
James Borger.
\newblock The basic geometry of {W}itt vectors, {I}: {T}he affine case.
\newblock {\em Algebra Number Theory}, 5(2):231--285, 2011.

\bibitem{Borger:BGWV-II}
James Borger.
\newblock The basic geometry of {W}itt vectors. {II}: {S}paces.
\newblock {\em Math. Ann.}, 351(4):877--933, 2011.

\bibitem{Buium:Diff-chars-over-p-adic}
Alexandru Buium.
\newblock Differential characters of abelian varieties over {$p$}-adic fields.
\newblock {\em Invent. Math.}, 122(2):309--340, 1995.

\bibitem{Buium:differential-modular-forms}
Alexandru Buium.
\newblock Differential modular forms.
\newblock {\em J. Reine Angew. Math.}, 520:95--167, 2000.

\bibitem{Erdogan:canonical-lifts}
Altan Erdo{\u{g}}an.
\newblock A universal formula for the j-invariant of the canonical lifting.
\newblock {\em Journal of Number Theory}, 150:26--40, 2015.

\bibitem{Finotti}
Lu{\'{\i}}s R.~A. Finotti.
\newblock Lifting the {$j$}-invariant: questions of {M}azur and {T}ate.
\newblock {\em J. Number Theory}, 130(3):620--638, 2010.

\bibitem{Thesis-Gurney}
Lance Gurney.
\newblock {\em Elliptic curves with complex multiplication and
  $\Lambda$-structures}.
\newblock PhD thesis, Australian National University, 2015.

\bibitem{Joyal:Witt}
Andr{\'e} Joyal.
\newblock {$\delta$}-anneaux et vecteurs de {W}itt.
\newblock {\em C. R. Math. Rep. Acad. Sci. Canada}, 7(3):177--182, 1985.

\bibitem{Katz:modular-forms}
Nicholas~M. Katz.
\newblock {$p$}-adic properties of modular schemes and modular forms.
\newblock In {\em Modular functions of one variable, {III} ({P}roc. {I}nternat.
  {S}ummer {S}chool, {U}niv. {A}ntwerp, {A}ntwerp, 1972)}, pages 69--190.
  Lecture Notes in Mathematics, Vol. 350. Springer, Berlin, 1973.

\bibitem{Katz-Mazur}
Nicholas~M. Katz and Barry Mazur.
\newblock {\em Arithmetic moduli of elliptic curves}, volume 108 of {\em Annals
  of Mathematics Studies}.
\newblock Princeton University Press, Princeton, NJ, 1985.

\bibitem{Langer-Zink:dRW}
Andreas Langer and Thomas Zink.
\newblock De {R}ham-{W}itt cohomology for a proper and smooth morphism.
\newblock {\em J. Inst. Math. Jussieu}, 3(2):231--314, 2004.

\bibitem{Lazard:formal-groups-book}
Michel Lazard.
\newblock {\em Commutative formal groups}.
\newblock Lecture Notes in Mathematics, Vol. 443. Springer-Verlag, Berlin-New
  York, 1975.

\bibitem{Messing:Barsotti-Tate}
William Messing.
\newblock {\em The crystals associated to {B}arsotti-{T}ate groups: with
  applications to abelian schemes}.
\newblock Lecture Notes in Mathematics, Vol. 264. Springer-Verlag, Berlin-New
  York, 1972.

\bibitem{stacks-project}
The {Stacks Project Authors}.
\newblock {\itshape {S}tacks {P}roject}.
\newblock http://stacks.math.columbia.edu, 2016.

\bibitem{van-der-Kallen:Descent}
Wilberd van~der Kallen.
\newblock Descent for the {$K$}-theory of polynomial rings.
\newblock {\em Math. Z.}, 191(3):405--415, 1986.

\end{thebibliography}
\bibliographystyle{plain}

\end{document}